# Classification of 7-dimensional subalgebras of the 8-dimensional Clifford algebra

## U. Shtukar

Author: Dr. Uladzimir Shtukar, Associate Professor, North Carolina Central University, USA.

e-mail: vshtukar@yahoo.com  Post address: 1906 Raj Drive, Durham, NC 27703.

**Abstract.** All 7-dimensional subalgebras of the 8-dimensional Clifford algebra over the field C of complex numbers are found. Canonical bases are used throughout the determination. It is found that the 8-dimensional Clifford algebra over C has exactly eight 7-dimensional subalgebras and each of these eight is over the complex numbers. It follows that the 8-dimensional Clifford algebra over C has no 7-dimensional subalgebra over the real numbers. Another consequence is that the eight 7-dimensional subalgebras of the 8-dimentional Clifford over C are in fact 7-dimensional subalgebras of all $2^n$-dimensional Clifford algebras over C for $n \geq 3$.

Let $V_n$ be the $n$−dimensional vector space with its standard basis $e_1, e_2, \ldots, e_n$. The Clifford algebra over $V_n$ is the associative algebra whose generators are $e_A = e_{i_1} e_{i_2} \ldots e_{i_k}$ for each possible index $A = (i_1 i_2 \ldots i_k)$ that is a subset of {1,2, …, n} with $i_1 < i_2 < \cdots < i_k$. It follows that (as a vector space) the dimension of the Clifford algebra over $V_n$ is $2^n$. The multiplicative identity element of the Clifford algebra is $e_\emptyset$. It is also denoted by $e_0$ and 1; and, it has the following property $e_0 e_i = e_i e_0 = e_i$, $e_0 e_0 = 1$. The vectors $e_1, e_2, \ldots, e_n$ are included in the generators of the Clifford algebra (via index sets A with k = 1) and the products of these particular vectors satisfy the following conditions

$$e_i e_j = -e_j e_i, (i \neq j), \text{ and } e_i^2 = -1, \text{ where } i, j = 1, 2, \ldots, n.$$

As mentioned, the multiplication operation of the Clifford algebra is associative, so $(ab)c = a(bc)$ for any elements $a, b, c$ of the Clifford algebra. An arbitrary element of the Clifford algebra over $V_n$ is

$$a = \sum_A a_A e_A,$$

using all indices $A = (i_1 i_2 \ldots i_k)$. We denote the Clifford algebra over $V_n$ by $g$(n,C), when the coefficients $a_A$ are complex numbers and by g(n,R), when the coefficients $a_A$ are real numbers. For more information about Clifford algebras, see, for example, [1].

Classification of subalgebras is a fundamental way to gain understanding of any algebra. What can we say about subalgebras of the Clifford algebra $g$(n,C)? First, it is easy to see that the Clifford algebra $g$(n,C) is a $2^n$−dimensional subalgebra of the $2^{n+k}$−dimensional Clifford algebra $g$(n+k,C) with $k \geq 0$. This leads to the more general and equally obvious statement.

**Lemma. If $h$ is a subalgebra of $2^n$−dimensional Clifford algebra $g$(n,C), then $h$ is also a subalgebra of $2^{n+k}$−dimensional Clifford algebra $g$(n+k,C), where $k \geq 0$.**

A lot of properties of Clifford algebras have been established, and some concern subalgebras. For example, equal rank subalgebras are discussed by E. Meinrenken in his book [2]. An answer to the following question would be real progress in the analysis of Clifford algebras: Does the Clifford algebra $g$(n,C) have $k$−dimensional subalgebras for $2^{n-1} < k < 2^n$? This article answers this question for the case $k=7$, $n=3$ by finding all 7-dimensional subalgebras of the 8−dimensional Clifford algebra $g(3,C)$.

For simplicity, the following notation will be used for the generators of $g(3,C)$:

$$e_0 = 1, e_1, e_2, e_3, i = e_1e_2, j = e_1e_3, k = e_2e_3, z = e_1e_2e_3.$$

Products of the generators of $g(3,C)$ are recorded in the following table:

| $\varnothing$ | 1 | $e_1$ | $e_2$ | $e_3$ | $i$ | $j$ | $k$ | $z$ |
|---|---|---|---|---|---|---|---|---|
| 1 | 1 | $e_1$ | $e_2$ | $e_3$ | $i$ | $j$ | $k$ | $z$ |
| $e_1$ | $e_1$ | $-1$ | $i$ | $j$ | $-e_2$ | $-e_3$ | $z$ | $-k$ |
| $e_2$ | $e_2$ | $-i$ | $-1$ | $k$ | $e_1$ | $-z$ | $-e_3$ | $j$ |
| $e_3$ | $e_3$ | $-j$ | $-k$ | $-1$ | $z$ | $e_1$ | $e_2$ | $-i$ |
| $i$ | $i$ | $e_2$ | $-e_1$ | $z$ | $-1$ | $k$ | $-j$ | $-e_3$ |
| $j$ | $j$ | $e_3$ | $-z$ | $-e_1$ | $-k$ | $-1$ | $i$ | $e_2$ |
| $k$ | $k$ | $z$ | $e_3$ | $-e_2$ | $j$ | $-i$ | $-1$ | $-e_1$ |
| $z$ | $z$ | $-k$ | $j$ | $-i$ | $-e_3$ | $e_2$ | $-e_1$ | 1 |

The following statement contains the total description of all 7-dimensional subalgebras of the 8-dimensional Clifford algebra $g(3,C)$.

**Theorem. The 8-dimensional Clifford algebra $g(3,C)$ has the following four 1-parameter sets of 7-dimensional subalgebras over the field of complex numbers $C$**

$$h_1 = Span\{1, e_1 + k, e_2 + ak, e_3 + \sqrt{-1-a^2}k, i - \sqrt{-1-a^2}k, j + ak, z\}, a \in C;$$

$$h_2 = Span\{1, e_1 + k, e_2 + ak, e_3 - \sqrt{-1-a^2}k, i + \sqrt{-1-a^2}k, j + ak, z\}, a \in C;$$

$$h_3 = Span\{1, e_1 - k, e_2 + ak, e_3 + \sqrt{-1-a^2}k, i + \sqrt{-1-a^2}k, j - ak, z\}, a \in C;$$

$$h_4 = Span\{1, e_1 - k, e_2 + ak, e_3 - \sqrt{-1-a^2}k, i - \sqrt{-1-a^2}k, j - ak, z\}, a \in C;$$

**and four isolated 7-dimensional subalgebras over the field of complex numbers $C$**

$$h_5 = Span\{1, e_1, e_2 + j, e_3 + \sqrt{-1}j, i + \sqrt{-1}j, k, z\},$$

$$h_6 = Span\{1, e_1, e_2 + j, e_3 - \sqrt{-1}j, i - \sqrt{-1}j, k, z\},$$

$$h_7 = Span\{1, e_1, e_2 - j, e_3 + \sqrt{-1}j, i - \sqrt{-1}j, k, z\},$$

$$h_8 = Span\{1, e_1, e_2 - j, e_3 - \sqrt{-1}j, i + \sqrt{-1}j, k, z\}.$$

**The Clifford algebra $g(3,C)$ has no 7-dimensional subalgebras over the field of real numbers.**

Proof. Since 7-dimensional subalgebras of $g(3,C)$ are 7-dimensional subspaces of $g(3,C)$, we first identify the 7-dimensional subspaces of $g(3,C)$ and then we identify which of these subspaces are subalgebras. To find 7-dimensional subspaces of the 8-dimensional $g(3,C)$, we use the canonical bases for $(n-1)$−dimensional subspaces of $n$−dimensional vector space that are detailed in [3]. Our case is $n=8$, and the corresponding canonical bases for this case are

(1) $a_1 = 1 + a_{18}z$, $a_2 = e_1 + a_{28}z$, $a_3 = e_2 + a_{38}z$, $a_4 = e_3 + a_{48}z$, $a_5 = i + a_{58}z$, $a_6 = j + a_{68}z$, $a_7 = k + a_{78}z$;

(2) $a_1 = 1 + a_{17}k$, $a_2 = e_1 + a_{27}k$, $a_3 = e_2 + a_{37}k$, $a_4 = e_3 + a_{47}k$, $a_5 = i + a_{57}k$, $a_6 = j + a_{67}k$, $a_7 = z$;

(3) $a_1 = 1 + a_{16}j$, $a_2 = e_1 + a_{26}j$, $a_3 = e_2 + a_{36}j$, $a_4 = e_3 + a_{46}j$, $a_5 = i + a_{56}j$, $a_6 = k$, $a_7 = z$;

(4) $a_1 = 1 + a_{15}i$, $a_2 = e_1 + a_{25}i$, $a_3 = e_2 + a_{35}i$, $a_4 = e_3 + a_{45}i$, $a_5 = j$, $a_6 = k$, $a_7 = z$;

(5) $a_1 = 1 + a_{14}e_3$, $a_2 = e_1 + a_{24}e_3$, $a_3 = e_2 + a_{34}e_3$, $a_4 = i$, $a_5 = j$, $a_6 = k$, $a_7 = z$;

(6) $a_1 = 1 + a_{13}e_2$, $a_2 = e_1 + a_{23}e_2$, $a_3 = e_3$, $a_4 = i$, $a_5 = j$, $a_6 = k$, $a_7 = z$;

(7) $a_1 = 1 + a_{12}e_1$, $a_2 = e_2$, $a_3 = e_3$, $a_4 = i$, $a_5 = j$, $a_6 = k$, $a_7 = z$;

(8) $a_1 = e_1$, $a_2 = e_2$, $a_3 = e_3$, $a_4 = i$, $a_5 = j$, $a_6 = k$, $a_7 = z$.

Every 7-dimensional subalgebra of $g(3,C)$ is a subspace of $g(3,C)$ and hence is associated with exactly one of these 8 canonical bases for some choice of the parameters $a_{ij}$. To determine which subspaces of $g(3,C)$ are subalgebras, we use the well-known fact, that a subspace $h$ of $g(3,C)$ is a subalgebra of $g(3,C)$ if and only if $xy \in h$ for any two elements $x \in h$ and $y \in h$. So, for each of the 8 canonical bases we find all products $a_i a_j$ for $i,j = 1,2,3,4,5,6,7$, and then determine whether or not all products $a_i a_j$ are in $h = Span\{a_1, a_2, a_3, a_4, a_5, a_6, a_7\}$. Specifically, for each of the 8 canonical bases we find all products $a_i a_j$ for $i,j = 1,2,3,4,5,6,7$, and then determine whether or not $x_1, x_2, x_3, x_4, x_5, x_6, x_7$ exist in C so that $a_i a_j = x_1 a_1 + x_2 a_2 + x_3 a_3 + x_4 a_4 + x_5 a_5 + x_6 a_6 + x_7 a_7$. The solutions to this equation are found by expressing each $a_i$ in terms of generators of $g(3,C)$ and then equating the coefficients from each side of the equation for corresponding generators. This technique either a) establishes conditions on the parameters $a_i a_j$ for a solution for the $x_i$ to exist or b) establishes that there is no solution. In case a) a subalgebra of $g(3,C)$ is associated with the canonical basis. In case b) no subalgebra of $g(3,C)$ is associated with the canonical basis.

**Basis (1).** Products $a_2 a_7$, $a_7 a_2$, $a_3 a_4$, $a_4 a_3$, $a_3 a_5$, $a_5 a_3$ play crucial roles for this basis. Start with $a_3 a_4$, and $a_4 a_3$. Using the table of products above, we have $a_3 a_4 = (e_2 + a_{38}z)(e_3 + a_{48}z) = k + a_{38}a_{48} + a_{48}j - a_{38}i$. To check whether or not this is in $h = Span\{a_1, a_2, a_3, a_4, a_5, a_6, a_7\}$, we look for $x_1, x_2, x_3, x_4, x_5, x_6, x_7$ so that $k + a_{38}a_{48} + a_{48}j - a_{38}i = x_1 a_1 + x_2 a_2 + x_3 a_3 + x_4 a_4 + x_5 a_5 + x_6 a_6 + x_7 a_7$. By expressing each $a_i$ in terms of generators of $g(3,C)$ and equating the coefficients from each side of the equation for corresponding generators, we find a solution exists ( $x_1 = a_{38}a_{48}$, $x_2 = 0$, $x_3 = 0$, $x_4 = 0$, $x_5 = -a_{38}$, $x_6 = a_{48}$, $x_7 = 1$ ), if $a_{18}a_{38}a_{48} - a_{38}a_{58} + a_{48}a_{68} + a_{78} = 0$. Similarly, the product $a_4 a_3$ requires the following condition for it to be in $h$: $a_{18}a_{38}a_{48} - a_{38}a_{58} + a_{48}a_{68} - a_{78} = 0$. From the last two equations, we must have $a_{78} = 0$ for both $a_3 a_4$, and $a_4 a_3$ to be in $h$.

Next, we have $a_3 a_5 = (e_2 + a_{38}z)(i + a_{58}z) = e_1 + a_{38}a_{58} - a_{38}e_3 + a_{58}j$. Using the same procedure that we used above, we find $a_3 a_5$ is in $h$, if $a_{18}a_{38}a_{58} + a_{28} - a_{38}a_{48} + a_{58}a_{68} = 0$. Similarly, the product $a_5 a_3$ requires the condition $a_{28} - a_{38}a_{48} + a_{58}a_{68} = 0$ for it to be in $h$. From the last two equations, we must have $a_{28} = 0$ for both $a_3 a_5$ and $a_5 a_3$ to be in $h$.

Now compute $a_2 a_7 = (e_1 + a_{28}z)(k + a_{78}z) = z + a_{28}a_{78} - a_{78}k - a_{28}e_1$. Using the same procedure, we find $a_2 a_7$ is in $h$, if $a_{18}a_{28}a_{78} - a_{28}^2 - a_{78}^2 = 1$. Product $a_7 a_2$ requires the same condition to be in $h$. If our previous requirements $a_{78} = 0$ and $a_{28} = 0$ are substituted into $a_{18}a_{28}a_{78} - a_{28}^2 - a_{78}^2 = 1$, the contradiction $0 = 1$ is produced. This means that it is impossible for all $a_i a_j$ to be in $h$. So, Basis (1) doesn't generate any subalgebra in algebra $g(3, C)$.

**Basis (2).** This basis consists of vectors $a_1 = 1 + a_{17}k$, $a_2 = e_1 + a_{27}k$, $a_3 = e_2 + a_{37}k$, $a_4 = e_3 + a_{47}k$, $a_5 = i + a_{57}k$, $a_6 = j + a_{67}k$, $a_7 = z$. Find all products $a_i a_j$ where $i, j = 1, 2, 3, 4, 5, 6, 7$, and find conditions so that all products $a_i a_j$ are in $h$ where $h = Span\{a_1, a_2, a_3, a_4, a_5, a_6, a_7\}$. Using the table of products, we have $a_1 a_7 = (1 + a_{17}k)z = z - a_{17}e_1$. To determine if this is in $h$, set it equal to $x_1 a_1 + x_2 a_2 + x_3 a_3 + x_4 a_4 + x_5 a_5 + x_6 a_6 + x_7 a_7$, express each $a_i$ in terms of the generators of $g(3, C)$, and solve for the $x_i$. This procedure yields a solution ($x_1 = x_3 = x_4 = x_5 = x_6 = 0$, $x_2 = -a_{17}$, $x_7 = 1$), if the condition $a_{17}a_{27} = 0$ is met. For the product $a_2 a_7$ we have $a_2 a_7 = (e_1 + a_{27}k)z = -k - a_{27}e_1$. By following the same technique, this is in $h$, if $a_{27}^2 = 1$. Under this condition, the previous condition ($a_{17}a_{27} = 0$) is equivalent to $a_{17} = 0$. Next, the product $a_3 a_7 = (e_2 + a_{37}k)z = j - a_{37}e_1$ is in $h$, if $-a_{27}a_{37} + a_{67} = 0$. The product $a_4 a_7 = (e_3 + a_{47}k)z = -i - a_{47}e_1$ is in $h$, if $a_{27}a_{47} + a_{57} = 0$. The product $a_5 a_7 = (i + a_{57}k)z = -e_3 - a_{57}e_1$ is in $h$, if $a_{27}a_{57} + a_{47} = 0$. This is not a new condition. Since $a_{27}^2 = 1$, it is equivalent to $a_{27}a_{47} + a_{57} = 0$ which is the condition for $a_4 a_7$ to be in $h$. The product $a_6 a_7 = (j + a_{67}k)z = e_2 - a_{67}e_1$ is in $h$, if $-a_{27}a_{67} + a_{37} = 0$. This condition is also not new. Since $a_{27}^2 = 1$, it is equivalent to $-a_{27}a_{37} + a_{67} = 0$ which is the condition for $a_3 a_7$ to be in $h$. The following initial list summarizes the conditions that the parameters $a_{17}, a_{27}, a_{37}, a_{47}, a_{57}, a_{67}$ must satisfy to insure that all products $a_i a_j$ considered thus far will be in $h = Span\{a_1, a_2, a_3, a_4, a_5, a_6, a_7\}$:

$$a_{17} = 0, \; a_{27}^2 = 1, \; a_{57} = -a_{27}a_{47}, \; \text{and} \; a_{67} = a_{27}a_{37}.$$

Next note that $a_1 a_1 = (1 + a_{17}k)(1 + a_{17}k) = 1 - a_{17}^2 + 2a_{17}k$. It turns out that this is in $h$, if $a_{17} = 0$ which is one of the conditions in our initial list. Similarly, we have $a_2 a_2 \in h$, $a_3 a_3 \in h$, $a_4 a_4 \in h$, $a_5 a_5 \in h$, $a_6 a_6 \in h$, $a_7 a_7 \in h$. Also, since $a_{17} = 0$, we have $a_1 = 1$ and, therefore, $a_1 a_2 = a_2 a_1 = a_2$, $a_1 a_3 = a_3 a_1 = a_3$, $a_1 a_4 = a_4 a_1 = a_4$, $a_1 a_5 = a_5 a_1 = a_5$, $a_1 a_6 = a_6 a_1 = a_6$. So, all of these products are in $h$ without requiring any additional conditions for parameters $a_{17}, a_{27}, a_{37}, a_{47}, a_{57}, a_{67}$.

Additionally, the products $a_7 a_1, a_7 a_2, a_7 a_3, a_7 a_4, a_7 a_5, a_7 a_6$ are in $h$ without requiring any additional conditions for parameters $a_{17}, a_{27}, a_{37}, a_{47}, a_{57}, a_{67}$. These products are equal to $a_1 a_7, a_2 a_7, a_3 a_7, a_4 a_7, a_5 a_7, a_6 a_7$ respectively and all of these products were used in constructing the initial list of conditions on the parameters.

Now consider the remaining products of basic elements for Basis (2). For each, we will use the technique employed above to find any condition that must be met for the product to belong to $h$.

The product $a_2 a_3 = (e_1 + a_{27}k)(e_2 + a_{37}k) = i + a_{37}z + a_{27}e_3 - a_{27}a_{37}$ is in $h$, if $a_{27}a_{47} + a_{57} = 0$. The same condition is required for $a_3 a_2$ to be in $h$. This condition is in the initial list.

Next, $a_2 a_4 = (e_1 + a_{27}k)(e_3 + a_{47}k) = j + a_{47}z - a_{27}e_2 - a_{27}a_{47}$ is in $h$, if $-a_{27}a_{37} + a_{67} = 0$. The same condition is required for $a_4 a_2$ to be in $h$. This condition is in the initial list.

Next, $a_2 a_5 = (e_1 + a_{27}k)(i + a_{57}k) = -e_2 + a_{57}z + a_{27}j - a_{27}a_{57}$ is in $h$, if $-a_{37} + a_{27}a_{67} = 0$. The same condition is required for $a_5 a_2$ to be in $h$. Since $a_{27}^2 = 1$, this condition is equivalent to $a_{67} = a_{27}a_{37}$ which is in the initial list.

Next, $a_2a_6 = (e_1 + a_{27}k)(j + a_{67}k) = -e_3 + a_{67}z - a_{27}i - a_{27}a_{67}$ is in $h$, if $-a_{47} - a_{27}a_{57} = 0$. The same condition is required for $a_6a_2$ to be in $h$. Since $a_{27}^2 = 1$, this condition is equivalent to $a_{57} = -a_{27}a_{47}$ which is in the initial list.

Next, $a_3a_4 = (e_2 + a_{37}k)(e_3 + a_{47}k) = k - a_{47}e_3 - a_{37}e_2 - a_{37}a_{47}$ is in $h$, if $-a_{37}a_{37} - a_{47}a_{47} = 1$. The same condition is required for $a_4a_3$ to be in $h$. This is a new condition.

Next, $a_3a_5 = (e_2 + a_{37}k)(i + a_{57}k) = e_1 - a_{57}e_3 + a_{37}j - a_{37}a_{57}$ is in h, if $a_{27} - a_{47}a_{57} + a_{37}a_{67} = 0$. The same condition is required for $a_5a_3$ to be in $h$. This is a new condition.

Next, $a_3a_6 = (e_2 + a_{37}k)(j + a_{67}k) = -z - a_{67}e_3 - a_{37}i - a_{37}a_{67}$ is in $h$, if $-a_{47}a_{67} - a_{37}a_{57} = 0$. The same condition is required for $a_6a_3$ to be in $h$. This is a new condition.

Next, $a_4a_5 = (e_3 + a_{47}k)(i + a_{57}k) = z - a_{57}e_2 + a_{47}j - a_{47}a_{57}$ is in $h$, if $a_{37}a_{57} + a_{47}a_{67} = 0$. The same condition is required for $a_5a_4$ to be in $h$. This condition is the same as the new condition found for $a_3a_6$.

Next, $a_4a_6 = (e_3 + a_{47}k)(j + a_{67}k) = e_1 + a_{67}e_2 - a_{47}i - a_{47}a_{67}$ is in $h$, if $a_{27} + a_{37}a_{67} - a_{47}a_{57} = 0$. The same condition is required for $a_6a_4$ to be in $h$. This condition is the same as the new condition found for $a_3a_5$.

Next, $a_5a_6 = (i + a_{57}k)(j + a_{67}k) = k - a_{67}j - a_{57}i - a_{57}a_{67}$ is in h, if $a_{57}a_{57} + a_{67}a_{67} = -1$. The same condition is required for $a_6a_5$ to be in $h$. This is equivalent to the new condition $-a_{37}a_{37} - a_{47}a_{47} = 1$ found for $a_3a_4$. To see this, note that the three conditions $a_{57} = -a_{27}a_{47}$, $a_{67} = a_{27}a_{37}$ and $a_{27}^2 = 1$ are in the initial list and together imply $a_{57}^2 = a_{47}^2$ and $a_{67}^2 = a_{37}^2$.

The following list summarizes the conditions that the parameters $a_{17}, a_{27}, a_{37}, a_{47}, a_{57}, a_{67}$ must satisfy to insure that all products $a_ia_j$ considered thus far will be in $h = Span\{a_1, a_2, a_3, a_4, a_5, a_6, a_7\}$ for $i, j = 1, 2, 3, 4, 5, 6, 7$:

$$a_{17} = 0, \ a_{27}^2 = 1, \ a_{57} = -a_{27}a_{47}, \ a_{67} = a_{27}a_{37}, \ -a_{37}a_{37} - a_{47}a_{47} = 1,$$

$$-a_{47}a_{67} - a_{37}a_{57} = 0, \quad a_{27} + a_{37}a_{67} - a_{47}a_{57} = 0.$$

If $a_{27} = 1$, then $a_{57} = -a_{47}$, $a_{67} = a_{37}$, and $a_{37}^2 + a_{47}^2 = -1$. The last equation has no solution in the field of real numbers but it has solutions in the field of complex numbers. For convenience, denote $a_{37} = a$. Then the parameter values $a_{17} = 0$, $a_{27} = 1$, $a_{37} = a$, $a_{47} = \pm\sqrt{-1-a^2}$, $a_{57} = \mp\sqrt{-1-a^2}$, and $a_{67} = a$ clearly satisfy the first five of the seven conditions that are listed above as sufficient to put all products $a_ia_j$ in $h = Span\{a_1, a_2, a_3, a_4, a_5, a_6, a_7\}$. It is easy to check that these parameter values also satisfy the last two conditions. Therefore, the following 1-parameter sets of subalgebras in algebra $g(3, C)$ are found:

$$h_1 = Span\{1, e_1 + k, e_2 + ak, e_3 + \sqrt{-1-a^2}k, i - \sqrt{-1-a^2}k, j + ak, z\}, a \in C;$$

$$h_2 = Span\{1, e_1 + k, e_2 + ak, e_3 - \sqrt{-1-a^2}k, i + \sqrt{-1-a^2}k, j + ak, z\}, a \in C.$$

If $a_{27} = -1$, then $a_{57} = a_{47}$, $a_{67} = -a_{37}$, $a_{37}^2 + a_{47}^2 = -1$. The last equation has no solution in the field of real numbers but it has solutions in the field of complex numbers. For convenience, denote $a_{37} = a$. Then the parameter values $a_{17} = 0$, $a_{27} = -1$, $a_{37} = a$, $a_{47} = \pm\sqrt{-1-a^2}$, $a_{57} = \pm\sqrt{-1-a^2}$ and $a_{67} = -a$ clearly satisfy the first five of the seven conditions that are listed above as sufficient to put all products $a_ia_j$ in $h = Span\{a_1, a_2, a_3, a_4, a_5, a_6, a_7\}$. It is easy to check that these parameter values also satisfy the last two conditions. Therefore, the following 1-parameter sets of subalgebras in algebra $g(3, C)$ are found:

$$h_3 = Span\{1, e_1 - k, e_2 + ak, e_3 + \sqrt{-1-a^2}k, i + \sqrt{-1-a^2}k, j - ak, z\}, a \in C;$$

$$h_4 = Span\{1, e_1 - k, e_2 + ak, e_3 - \sqrt{-1-a^2}k, i - \sqrt{-1-a^2}k, j - ak, z\}, a \in C.$$

**Basis (3).** This basis consists of vectors $a_1 = 1 + a_{16}j$, $a_2 = e_1 + a_{26}j$, $a_3 = e_2 + a_{36}j$, $a_4 = e_3 + a_{46}j$, $a_5 = i + a_{56}j$, $a_6 = k$, $a_7 = z$. Consider subspace $h = Span\{a_1, a_2, a_3, a_4, a_5, a_6, a_7\}$. Find all products $a_i a_j$ where $i, j = 1, 2, 3, 4, 5, 6, 7$, and find conditions so that all products $a_i a_j$ are in $h$. Using the table of products, we have $a_1 a_7 = (1 + a_{16}j)z = z + a_{16}e_2$. This turns out to be in $h$, if $a_{16}a_{36} = 0$. For the product $a_2 a_7$ we have $a_2 a_7 = (e_1 + a_{26}j)z = -k + a_{26}e_2$. This is in $h$, if $a_{26}a_{36} = 0$. Next, compute the product $a_3 a_7 = (e_2 + a_{36}j)z = j + a_{36}e_2$. This is in $h$, if $a_{36}^2 = 1$. Under this new condition, the two previous conditions ($a_{16}a_{36} = 0$ and $a_{26}a_{36} = 0$) are equivalent to requiring $a_{16} = 0$ and $a_{26} = 0$. Now compute product $a_4 a_7 = (e_3 + a_{46}j)z = -i + a_{46}e_2$. This is in $h$, if $a_{36}a_{46} - a_{56} = 0$. Next, compute $a_5 a_7 = (i + a_{56}j)z = -e_3 + a_{56}e_2$. This is in $h$, if $a_{36}a_{56} - a_{46} = 0$. This condition is not new. Since $a_{36}^2 = 1$, it is equivalent to $a_{36}a_{46} - a_{56} = 0$ which is the condition for $a_4 a_7$ to be in $h$. Compute product $a_6 a_7 = kz = -e_1$. This is in $h$, if $a_{26} = 0$. This is also a condition we have already found. The following initial list summarizes the conditions that the parameters $a_{16}, a_{26}, a_{36}, a_{46}, a_{56}$ must satisfy to insure that all products $a_i a_j$ considered thus far will be in $h = Span\{a_1, a_2, a_3, a_4, a_5, a_6, a_7\}$:

$$a_{16} = 0, \; a_{26} = 0, \; a_{36}^2 = 1, \text{ and } a_{56} = a_{36}a_{46}.$$

Next note that $a_1 a_1 = (1 + a_{16}j)(1 + a_{16}j) = 1 - a_{16}^2 + 2a_{16}j$. It turns out that this is in $h$, if $a_{16} = 0$ which is one of the conditions listed above. Similarly, we have $a_2 a_2 \in h$, $a_3 a_3 \in h$, $a_4 a_4 \in h$, $a_5 a_5 \in h$, $a_6 a_6 \in h$, $a_7 a_7 \in h$. Also, since $a_{16} = 0$, we have $a_1 = 1$ and, therefore, $a_1 a_2 = a_2 a_1 = a_2$, $a_1 a_3 = a_3 a_1 = a_3$, $a_1 a_4 = a_4 a_1 = a_4$, $a_1 a_5 = a_5 a_1 = a_5$, $a_1 a_6 = a_6 a_1 = a_6$. So, all of these products are in $h$ without requiring any additional conditions for parameters $a_{16}$, $a_{26}$, $a_{36}$, $a_{46}$, $a_{56}$.

Additionally, the products $a_7 a_1$, $a_7 a_2$, $a_7 a_3$, $a_7 a_4$, $a_7 a_5$, $a_7 a_6$ are in $h$ without requiring any additional conditions for parameters $a_{16}$, $a_{26}$, $a_{36}$, $a_{46}$, $a_{56}$. These products are equal to $a_1 a_7$, $a_2 a_7$, $a_3 a_7$, $a_4 a_7$, $a_5 a_7$, $a_6 a_7$ respectively and all of these products were used in constructing the initial list of conditions on the parameters.

Now consider the remaining products of basic elements for Basis (3). For each, we will use the technique employed above to find any condition that must be met for the product to belong to $h$.

The product $a_2 a_3 = (e_1 + a_{26}j)(e_2 + a_{36}j) = i - a_{36}e_3 - a_{26}z - a_{26}a_{36}$ is in $h$, if $-a_{36}a_{46} + a_{56} = 0$. The same condition is required for $a_3 a_2$ to be in $h$. This condition is in the initial list above.

Next, $a_2 a_4 = (e_1 + a_{26}j)(e_3 + a_{46}j) = j - a_{46}e_3 - a_{26}e_1 - a_{26}a_{46}$ is in $h$, if $-a_{26}a_{26} - a_{46}a_{46} = 1$. The same condition is required for $a_4 a_2$ to be in $h$. Since $a_{26} = 0$, this is equivalent to $a_{46}a_{46} = -1$. This is a new condition.

Next, $a_2 a_5 = (e_1 + a_{26}j)(i + a_{56}j) = -e_2 - a_{56}e_3 - a_{26}k - a_{26}a_{56}$ is in $h$, if $-a_{36} - a_{46}a_{56} = 0$. The same condition is required for $a_5 a_2$ to be in $h$. Since $a_{46}a_{46} = -1$ (the new condition for $a_2 a_4$), this condition is equivalent to $a_{56} = a_{36}a_{46}$ which is in the initial list.

Next, $a_2 a_6 = (e_1 + a_{26}j)k = z + a_{26}i$ is in $h$, if $a_{26}a_{56} = 0$. The same condition is required for $a_6 a_2$ to be in $h$. This condition is a consequence of $a_{26} = 0$ which is in the initial list.

Next, $a_3 a_4 = (e_2 + a_{36}j)(e_3 + a_{46}j) = k - a_{46}z - a_{36}e_1 - a_{36}a_{46}$ is in $h$, if $a_{26}a_{36} = 0$. The same condition is required for $a_4 a_3$ to be in $h$. This condition is a consequence of $a_{26} = 0$ which is in the initial list.

Next, $a_3 a_5 = (e_2 + a_{36}j)(i + a_{56}j) = e_1 - a_{56}z - a_{36}k - a_{36}a_{56}$ is in $h$, if $a_{26} = 0$. Product $a_5 a_3$ requires the same condition to be in $h$. This condition is in the initial list.

Next, $a_3 a_6 = (e_2 + a_{36}j)k = -e_3 + a_{36}i$ is in $h$, if $-a_{46} + a_{36}a_{56} = 0$. Product $a_6 a_3$ requires the same condition to be in $h$. Since $a_{36}^2 = 1$, this condition is equivalent to $a_{56} = a_{36}a_{46}$ which is in the initial list.

Next, $a_4 a_5 = (e_3 + a_{46}j)(i + a_{56}j) = z + a_{56}e_1 - a_{46}k - a_{46}a_{56}$ is in $h$, if $a_{26}a_{56} = 0$. Similarly, product $a_5 a_4$ requires the same condition to be in $h$. This condition is a consequence of $a_{26} = 0$ which is in the initial list.

Next, $a_4 a_6 = (e_3 + a_{46}j)k = e_2 + a_{46}i$ is in $h$, if $a_{36} + a_{46}a_{56} = 0$. Product $a_6 a_4$ requires the same condition to be in $h$. Since $a_{46}a_{46} = -1$ (the new condition for $a_2 a_4$), this condition is equivalent to $a_{56} = a_{36}a_{46}$ which is in the initial list.

Next, $a_5 a_6 = (i + a_{56}j)k = -j + a_{56}i$ is in $h$, if $a_{56}a_{56} = -1$. Product $a_6 a_5$ requires the same condition to be in $h$. This is a new condition.

The following list summarizes the conditions that the parameters $a_{16}, a_{26}, a_{36}, a_{46}, a_{56}$ must satisfy to insure that all products $a_i a_j$ considered thus far will be in $h = Span\{a_1, a_2, a_3, a_4, a_5, a_6, a_7\}$ for $i, j = 1, 2, 3, 4, 5, 6, 7$:

$$a_{16} = 0, \ a_{26} = 0, \ a_{36}^2 = 1, \ a_{56} = a_{36}a_{46}, \ a_{46}a_{46} = -1, \ a_{56}a_{56} = -1.$$

If $a_{36} = 1$, then the parameter values $a_{16} = 0$, $a_{26} = 0$, $a_{36} = 1$, $a_{46} = \pm\sqrt{-1}$, and $a_{56} = a_{46} = \pm\sqrt{-1}$ clearly satisfy the six conditions that are listed above as sufficient to put all products $a_i a_j$ in $h = Span\{a_1, a_2, a_3, a_4, a_5, a_6, a_7\}$. If $a_{36} = -1$, then the parameter values $a_{16} = 0$, $a_{26} = 0$, $a_{36} = -1$, $a_{46} = \pm\sqrt{-1}$, and $a_{56} = -a_{46} = \mp\sqrt{-1}$ clearly satisfy the six conditions that are listed above as sufficient to put all products $a_i a_j$ in $h = Span\{a_1, a_2, a_3, a_4, a_5, a_6, a_7\}$. Therefore, the system has four solutions, and four subalgebras of algebra $g(3, C)$ are found:

$$h_5 = Span\{1, e_1, e_2 + j, e_3 + \sqrt{-1}j, i + \sqrt{-1}j, k, z\},$$

$$h_6 = Span\{1, e_1, e_2 + j, e_3 - \sqrt{-1}j, i - \sqrt{-1}j, k, z\},$$

$$h_7 = Span\{1, e_1, e_2 - j, e_3 + \sqrt{-1}j, i - \sqrt{-1}j, k, z\},$$

$$h_8 = Span\{1, e_1, e_2 - j, e_3 - \sqrt{-1}j, i + \sqrt{-1}j, k, z\}.$$

**Basis (4).** This basis consists of vectors $a_1 = 1 + a_{15}i$, $a_2 = e_1 + a_{25}i$, $a_3 = e_2 + a_{35}i$, $a_4 = e_3 + a_{45}i$, $a_5 = j$, $a_6 = k$, $a_7 = z$. Consider the subspace $h = Span\{a_1, a_2, a_3, a_4, a_5, a_6, a_7\}$ and the specific products $a_5 a_6$ and $a_6 a_5$. We have $a_5 a_6 = jk = i$, and $a_6 a_5 = kj = -i$. But, it is easy to check that vectors $i$ and $-i$ don't belong to the subspace $h$. Therefore, $h$ is not a subalgebra of algebra $g(3, C)$.

**Basis (5).** This basis consists of vectors $a_1 = 1 + a_{14}e_3$, $a_2 = e_1 + a_{24}e_3$, $a_3 = e_2 + a_{34}e_3$, $a_4 = i$, $a_5 = j$, $a_6 = k$, $a_7 = z$. Consider the subspace $h = Span\{a_1, a_2, a_3, a_4, a_5, a_6, a_7\}$ and the specific products $a_4a_7$ and $a_7a_4$. We have $a_4a_7 = iz = -e_3$, and $a_7a_4 = zi = -e_3$. But, it is easy to check that vector $e_3$ doesn't belong to the subspace $h$. Therefore, $h$ is not a subalgebra of algebra $g(3,C)$.

**Bases (6), (7), (8).** These bases don't generate any subalgebra of algebra $g(3,C)$. These three cases are very similar to the cases of Bases (4) and (5); so, the details are omitted.

The proof is finished.

**Remark.** According to Lemma, we can say that the eight 7-dimensional subalgebras of $g(3,C)$ found here are 7-dimensional subalgebras of each $2^n$ −dimensional Clifford algebra $g(n,C)$ for $n \geq 3$.

This article begins the classification of small dimensional subalgebras for Clifford algebras. The results are obtained by using canonical bases.

Author thanks Professor Russell Gosnell for his useful advise.

**References**

[1] J. Chisholm, A. Common. "Clifford Algebras and Their Applications in Mathematical Physics". Reidel, Dordrecht, 1986.

[2] E. Meinrenken. "Clifford Algebras and Lie Theory". Springer, Berlin – Heidelberg, 2013.

[3] U. Shtukar. "Classification of Canonical Bases for $(n-1)$ −Dimensional Subspaces of $n$ −Dimensional Vector Space". Journal of Generalized Lie Theory and Applications. Vol. 10, Issue 1, 2016, 241-245.